\documentclass[12pt]{amsart}
\usepackage{pythonhighlight} 
\usepackage{graphicx, breqn}
\usepackage{pythonhighlight} 
\usepackage{setspace}
\usepackage{subcaption}
\usepackage{amssymb}
\usepackage{multirow}
\usepackage{longtable}

\usepackage{hyperref}

\usepackage{amsthm}
\usepackage[mathscr]{eucal}
\usepackage{float, url}
\usepackage{amsmath, amsthm, tikz, fullpage, longtable,ulem,amssymb}
\usepackage{xpatch}
\usetikzlibrary{arrows.meta, shadows, fadings, shapes.arrows, positioning}
\tikzset{My Arrow Style/.style={single arrow, draw, text width=0.75cm}}
\makeatletter

\newcounter{proofpart}
\xpretocmd{\proof}{\setcounter{proofpart}{0}}{}{}
\newcommand{\proofpart}[1]{%
  \par
  \addvspace{\medskipamount}%
  \stepcounter{proofpart}%
  \noindent\emph{Part \theproofpart: #1}\par\nobreak\smallskip
  \@afterheading
}

\usepackage{color}
\usepackage{listings}
\usepackage{setspace}

\lstdefinelanguage{Python}{
numbers=left,
numberstyle=\footnotesize,
numbersep=1em,
xleftmargin=1em,
framextopmargin=2em,
framexbottommargin=2em,
showspaces=false,
showtabs=false,
showstringspaces=false,
frame=l,
tabsize=4,
basicstyle=\ttfamily\small\setstretch{1},
backgroundcolor=\color{Background},
commentstyle=\color{Comments}\slshape,
stringstyle=\color{Strings},
morecomment=[s][\color{Strings}]{"""}{"""},
morecomment=[s][\color{Strings}]{'''}{'''},
morekeywords={import,from,class,def,for,while,if,is,in,elif,else,not,and,or,print,break,continue,return,True,False,None,access,as,,del,except,exec,finally,global,import,lambda,pass,print,raise,try,assert},
keywordstyle={\color{Keywords}\bfseries},
morekeywords={[2]@invariant,pylab,numpy,np,scipy},
keywordstyle={[2]\color{Decorators}\slshape},
emph={self},
emphstyle={\color{self}\slshape},
}

\makeatother

\theoremstyle{definition}

\newtheorem{theorem}{Theorem}
\newtheorem{lemma}[theorem]{Lemma}
\newtheorem{corr}[theorem]{Corollary}
\newtheorem{prop}[theorem]{Proposition}

\newtheorem{definition}{Definition}

\newcommand{\rc}[1]{\textcolor{red}{#1}}

\newcommand{\floor}[1]{\left\lfloor #1 \right\rfloor}

\newcommand{\ceil}[1]{\left\lceil #1 \right\rceil}

\newcommand{\hidethis}[1]{}

\newcommand{\wfG}{\omega(f_G)}
\newcommand{\wf}{\omega(f)}

\newcommand{\kaRD}{$k$aRD }
\newcommand{\taRD}{2aRD }
\newcommand{\oaRD}{1aRD }

\newcommand{\pkaRD}{$k$aRD}
\newcommand{\ptaRD}{2aRD}
\newcommand{\poaRD}{1aRD}

\newcommand{\gtaRD}[1]{\gamma_{{2\scriptstyle{\text{aRD}}}}(#1)}
\newcommand{\gkaRD}[1]{\gamma_{{k\scriptstyle{\text{aRD}}}}(#1)}

\newcommand{\GtaRD}[1]{\Gamma_{{2\scriptstyle{\text{aRD}}}}(#1)}

\newcommand{\gR}[1]{\gamma_R(#1)}
\newcommand{\lG}{\ell_{G}}
\newcommand{\lGG}{\ell_{G'}} 

\newcommand{\w}[1]{w(#1)}

\begin{document}
\thispagestyle{empty}
\title{On the $k$-attack Roman Dominating Number of a Graph}

\author{Garrison Koch}
\address{Rochester Institute of Technology}
\curraddr{1 Lomb Memorial Dr, Rochester, NY 14623}
\email{glk5534@g.rit.edu}

\author{Nathan Shank}
\address{Moravian University}
\curraddr{1200 Main Street, Bethlehem, PA 18018}
\email{shankn@moravian.edu}

\begin{abstract}
\noindent  Given a graph $G=(V,E)$, the dominating number of a graph is the minimum size of a vertex set, $V' \subseteq V$, so that every vertex in the graph is either in $V'$ or is adjacent to a vertex in $V'$. A Roman Dominating function of $G$ is defined as $f:V \rightarrow \{0,1,2\}$ such that every vertex with a label of 0 in $G$ is adjacent to a vertex with a label of 2. The Roman Dominating number of a graph is the minimum total weight over all possible Roman Dominating functions. We consider the $k$-attack Roman Domination, particularly focusing on 2-attack Roman Domination. A Roman Dominating function of $G$ is a $k$-attack Roman Dominating function of $G$ if for all $j\leq k$, any subset $S$ of $j$ vertices all with label 0 must have at least $j$ vertices with label 2 in the open neighborhood of $S$. The $k$-attack Roman Dominating number of $G, \gkaRD{G}$, is the minimum total weight over all possible $k$-attack Roman Dominating functions. We find $\gtaRD{G}$ for particular graph class, discuss properties of $k$-attack Roman Domination, and make several connections with other domination ideas. 
\end{abstract}

\maketitle

\renewcommand{\baselinestretch}{1.2}

\section{Introduction}

We consider only finite simple graphs $G=(V,E)$ where $V$ is a vertex set and $E$ is the edge set. Let $N(v)=\{u \in V: uv \in E\}$ denote the open neighborhood of the vertex $v$ and for any set of vertices, $S$, let $N(S) = \bigcup_{v \in S} N(v)$ be the open neighborhood of the set $S$. Let $N[v]= N(v) \cup \{v\}$ and $N[S]=N(S)\cup S$ denote the closed neighborhoods of $v$ and $S$ respectively.  For other common graph theory notation and terminology used throughout, we refer the reader to \cite{Buckley}. 

A set $V' \subseteq V$ is a \textit{dominating set} if every vertex $v \in V$ is in $V'$ or adjacent to a vertex in $V'$.  The minimum size of all dominating sets is called the \textit{dominating number} of the graph and is denoted $\gamma(G)$. Any dominating set, $V'$ with $|V'| = \gamma(G)$ is called a \textit{minimum dominating set}. The study of domination sets began with Berge \cite{Berge1958} in 1958 and Ore \cite{Ore1962} in 1962.  The topic picked up momentum in the late 1970's by the work of Cockayne and several co-authors (\cite{Cockayne1978}, \cite{Cockayne1975}, \cite{Cockayne1976}, and \cite{Cockayne1977}). A generalization of domination, \textit{Roman Domination} (RD), was introduced by Revelle and Rosing \cite{Revelle2000} for a perfect defense strategy modeled as a set-covering problem. Cockayne et al. \cite{Cockayne} define a \textit{RD function} on $G$ as a labeling $f_G:V \rightarrow \{0,1,2\}$ such that every vertex with a label of 0 in $G$ is adjacent to a vertex with a label of 2. The \textit{weight} of a RD function is $ \wfG = \sum_{v \in V(G)} f_G(v)$.  The \textit{RD number} of a graph $G$, denoted $\gR{G}$, is the minimum weight of all possible RD functions on $G$.  If $\wfG = \gR{G}$ then $f$ is a minimum weight RD function. For a survey of RD, see for example \cite{Chellali2020}, \cite{Chellali2021}, \cite{Favaron2009}, \cite{DominationBook2020}, and \cite{Haynes}. 

Cockayne et al. \cite{Cockayne} proved several bounds and properties for the RD functions.  An easy-to-see bound, $\gamma(G) \leq \gR{G} \leq 2\gamma(G)$, can be obtained by labeling all vertices in the minimum dominating set of $G$ with a label of 2. This bound does not hold for our new $k$-attack RD number, as we will soon see. Cockayne et al. also find general RD number for simple graphs like paths, trees, and cycles as we will do.  In addition they found a lower bound based on the order of the graph and the maximum degree.  

Roman Domination differs from domination in that we are allowed to store additional resources, or armies, at certain vertices. In a defense setting for RD, any vertex with a label 0 is unable to protect itself, any vertex with a label 1 is able to defend itself, and any vertex with a label 2 or greater can defend itself and help defend its neighbors. So RD changes how resources are distributed on a graph.  

RD variants we consider functions $f$ defined on $G$ which satisfy some conditions. Then we want to find the minimum weight over all such functions.  

Beeler et al. \cite{DRD} analyzed Double RD which consider functions $f_G:V \rightarrow \{0,1,2,3\}$ on graphs $G$ such that if $f(v) = 0$, then $v$ must be adjacent to at least two vertices with label 2 or one vertex with label 3 and if $f(v) =1$, then $v$ must be adjacent to at least one neighbor with label 2 or larger. The idea here being we want to defend against an attack of size two, possibly at the same vertex. Double RD was generalized to Triple RD (\cite{Abdollahzadeh2021}) and $[k]$-RD (\cite{Khalili2023},\cite{Valenzuela-Tripodoro2024}) using functions $f:V \rightarrow \{0, 1, \ldots, k+1\}$ and requiring that for any vertex $v$ with $f(v) <k$ then $f(v) + \sum_{u \in N^+(v)}(f(v)-1) \geq k$ where $N^+(v)$ are the neighbors, $u$, of $v$ with $f(u)\geq 1$. Thus an attack on $k$, not necessarily district, vertices can be defended.  

Roman $k$-Dominating was defined by K\"ammerlineg and Volkmann \cite{Volkmann2009} to defend against multiple attacks; however, our vertices are resource limited. A function $f:V\rightarrow \{0, 1, 2\}$ is a Roman $k$-Dominating function if for every vertex $v$ with $f(v) = 0$ there are at least $k$ vertices adjacent to $v$ with weight 2.  In this variant, any attack on $k$ vertices, not necessarily distinct, can be defended, however the vertices are resource limited. The minimum weight over all functions $f$ is the Roman $k$-Domination number and is denoted $\gamma_{kR}(G)$.

Strong RD was introduced by \'Alvarez-Ruiz et al. \cite{ALVAREZ2017}.  Strong RD considers a function $f:V \rightarrow \{0, 1, \ldots, \ceil{\frac{\Delta}{2}} + 1 \}$ where $\Delta$ is the maximum degree of $G$. If $V_0 = \{v \in V: f(v) = 0\}$ then $f$ is a Strong RD function if every $v \in V_0$ has a neighbor $w$ so that $f(w) \geq 2$ and $f(w)\geq 1 + \ceil{\frac{1}{2}|N(w) \cap V_0|}$. The idea here is that any vertex with a weight of at least 2 must be able to help defend at least half of its neighbors.  

Liu et. al. \cite{LIU2020} extended the idea of Strong RD to $k$-Strong RD.  Define a function $f:V \rightarrow \{0, 1, \ldots 1+\min\{\Delta, k\}\}$ so that zero-labeled vertices are undefended, 1-labeled vertices can only defend themselves, and any vertex with $f(u)\geq 2$ has the ability to transfer resources to at most $f(u)-1$ of its neighbors to defend them. A function $f$ is a $k$-Strong RD function if any attack on $k$ distinct vertices is able to be defended. The minimum weight over all functions $f$ is the $k$-Strong RD number and is denoted $\gamma_{k-SRD}(G)$.

Several other variants of RD have been explored. These variants often change the restrictions on the function $f$, if attacks must occur at distinct vertices, and/or the conditions on vertex weights relative to their neighbors. Some variants include $k$-Roman Domination (\cite{Henning2003}), perfect Roman Domination \cite{Klostermeyer}, independent Roman Domination \cite{Adabi2012}, maximal Roman Domination \cite{Ahanger2017}, mixed Roman Domination \cite{Ahanger2017b}, Roman \{2\}-Domination \cite{Chellali2015}, Roman $\{k\}$-Domination (\cite{Wang2021}), total Roman \{2\}-Domination \cite{Ahangar2022}, and weak Roman Domination \cite{HenningNew2003}, although this is not meant to be an exhaustive list.

Our new \textit{$k$-attack RD} will combine Roman $k$-dominating and $k$-Strong RD by defending against attacks on at most $k$ vertices; however, the vertices will be resource limited, meaning our function $f:V \rightarrow \{0, 1, 2\}$. 

The tactical idea we are modeling with the $k$-attack RD (\pkaRD) involves having up to $k$ distinct vertices attacked, however, we are only allow a limited number of resources at each vertex,  $(0, 1, \text{ or } 2)$. Thus, 1-attack RD (\poaRD) means that our enemy has one army and will attack one vertex, which is equivalent to the original RD. However, if we consider 2-attack RD (\ptaRD) we must fortify against two armies simultaneously attacking at two different vertices. In Figure \ref{fig:ALLRD} we see the same graph with different minimal labelings depending on the RD variant.

\begin{figure}[H]
\centering
\begin{minipage}{0.24\textwidth}
\centering
\tikzstyle{vertex}=[circle,fill=black,inner sep=4.5pt]
\tikzstyle{white}=[circle,fill=white,inner sep=1pt]
\tikzstyle{green}=[circle,fill=green,inner sep=4.5pt]
\begin{tikzpicture}
 \node [draw, white] (1) at (0,0) {0};  
 \node [draw, white] (2) at (0,1.5) {0};
 \node [draw, white] (3) at (.75,.75) {2};
 \node [draw, white] (4) at (1.5,1.5) {0};  
 \node [draw, white] (5) at (1.5,0) {0};
 \node [draw, white] (6) at (2.25,.75) {1};
 \node [draw, white] (7) at (1.5,.75) {0};
 \draw (1) -- (3); 
 \draw (2) -- (3); 
 \draw (3) -- (4); 
 \draw (3) -- (5); 
 \draw (4) -- (6); 
 \draw (5) -- (6); 
 \draw (7) -- (6);
 \draw (7) -- (3);
 \node[anchor=north] at (current bounding box.south){$\gamma_R(G) = 3$};
\end{tikzpicture}
\end{minipage}
\begin{minipage}{0.24\textwidth}
\centering
\tikzstyle{vertex}=[circle,fill=black,inner sep=4.5pt]
\tikzstyle{white}=[circle,fill=white,inner sep=1pt]
\tikzstyle{green}=[circle,fill=green,inner sep=4.5pt]
\begin{tikzpicture}
 \node [draw, white] (1) at (0,0) {0};  
 \node [draw, white] (2) at (0,1.5) {0};
 \node [draw, white] (3) at (.75,.75) {3};
 \node [draw, white] (4) at (1.5,1.5) {0};  
 \node [draw, white] (5) at (1.5,0) {0};
 \node [draw, white] (6) at (2.25,.75) {1};
 \node [draw, white] (7) at (1.5,.75) {0};
 \draw (1) -- (3); 
 \draw (2) -- (3); 
 \draw (3) -- (4); 
 \draw (3) -- (5); 
 \draw (4) -- (6); 
 \draw (5) -- (6);
 \draw (7) -- (6);
 \draw (7) -- (3);
 \node[anchor=north] at (current bounding box.south){$\gamma_{2-SRD}(G)=4$};
\end{tikzpicture}
\end{minipage}
\begin{minipage}{0.24\textwidth}
\centering

\tikzstyle{vertex}=[circle,fill=black,inner sep=4.5pt]
\tikzstyle{white}=[circle,fill=white,inner sep=1pt]
\tikzstyle{green}=[circle,fill=green,inner sep=4.5pt]
\begin{tikzpicture}
 \node [draw, white] (1) at (0,0) {1};  
 \node [draw, white] (2) at (0,1.5) {0};
 \node [draw, white] (3) at (.75,.75) {2};
 \node [draw, white] (4) at (1.5,1.5) {0};  
 \node [draw, white] (5) at (1.5,0) {0};
 \node [draw, white] (6) at (2.25,.75) {2};
 \node [draw, white] (7) at (1.5,.75) {0};
 \draw (1) -- (3); 
 \draw (2) -- (3); 
 \draw (3) -- (4); 
 \draw (3) -- (5); 
 \draw (4) -- (6); 
 \draw (5) -- (6); 
 \draw (7) -- (6);
 \draw (7) -- (3);
 \node[anchor=north] at (current bounding box.south){$\gtaRD{G}=5$};
\end{tikzpicture}
\end{minipage}
\begin{minipage}{0.24\textwidth}
\centering
\tikzstyle{vertex}=[circle,fill=black,inner sep=4.5pt]
\tikzstyle{white}=[circle,fill=white,inner sep=1pt]
\tikzstyle{green}=[circle,fill=green,inner sep=4.5pt]
\begin{tikzpicture}
 \node [draw, white] (1) at (0,0) {1};  
 \node [draw, white] (2) at (0,1.5) {1};
 \node [draw, white] (3) at (.75,.75) {2};
 \node [draw, white] (4) at (1.5,1.5) {0};  
 \node [draw, white] (5) at (1.5,0) {0};
 \node [draw, white] (6) at (2.25,.75) {2};
 \node [draw, white] (7) at (1.5,.75) {0};
 \draw (1) -- (3); 
 \draw (2) -- (3); 
 \draw (3) -- (4); 
 \draw (3) -- (5); 
 \draw (4) -- (6); 
 \draw (5) -- (6); 
 \draw (7) -- (6);
 \draw (7) -- (3);
 \node[anchor=north] at (current bounding box.south){$\gamma_{2R}(G)=6$};
\end{tikzpicture}
\end{minipage}
\caption{Valid minimum labelings for $G$ for RD (3), $2$-Strong RD (4), 2-attack RD (5), and Roman 2-Dominating (6).}
\label{fig:ALLRD}
\end{figure}

Applications of variants of RD could include planning for disaster relief, supply chain disruptions (\cite{LIU2020}), planning for medical intervention facilities and supplies, and complex coordinated terrorist attacks (CCTA). 

\section{Definitions}

The \textit{private neighbors} of a vertex $v$ relative to a set $S$ is $pn(v, S) = \{u \in V : N(u) \cap S = \{v\}\}$. The \textit{exterior private neighbors} of $v$ relative to $S$ is $epn(v, S) = pn(v, S) \cap (V - S)$. A vertex $v$ is an \textit{exterior private neighbor} of a set $S$ if $v \notin S$ and $v$ are adjacent to exactly one vertex in $S$. The set of exterior private neighbors of a set $S$ will be denoted by $epn(S)$. 

We now define the \textit{public neighbors} of a set $S \subseteq V$ as $\overline{p}n(S) = \{u \in V : |N(u) \cap S| \geq 2\}$, that is, $u$ is in the \textit{public neighborhood} of $S$ if $u$ is adjacent to at least two vertices in $S$. The \textit{exterior public neighbors} of a set $S$ as $\overline{ep}n(S) = \overline{p}n(v, S) \cap (V - S)$.

We will consider a \textit{weight function} as a vertex labeling $f_G:V(G) \rightarrow \{0, 1, 2\}$. Given a weight function $f_G$, denote the set of vertices in $G$ labeled 0 as $V_0(f_G)$, the set of vertices in $G$ labeled 1 as $V_1(f_G)$, and the set of vertices in $G$ labeled 2 as $V_2(f_G)$. When the underlying graph $G$ and weight function $f$ are well understood, we will write $f, V_0, V_1,$ and $V_2$.  Note that $f$ induces a partition $(V_0, V_1, V_2)$ of $V$.  The \textit{weight of the function $f_G$} is $\wfG = \sum_{v \in V(G)} f_G(v)$. 

Exterior public and exterior private neighborhoods will be  necessary as we will distinguish between vertices in $V_0$ which are (exterior) private neighbors of $V_2$ and vertices in $V_0$ that are (exterior) public neighbors of $V_2$.

\begin{definition}\label{defnkRD}
A \textit{$k$-attack RD function} (\kaRD function) of $G$ is a weight function $f_G:V\rightarrow\{0,1,2\}$ such that for any set $S \subseteq V_0$ with $|S| \leq k$, there is a set $T \subseteq V_2 \cap N(S)$ with $|T| \geq |S|$. 
\end{definition}

Intuitively, we can consider an enemy with at most $k$ armies, and our resources are distributed based on a weight function $f$. Vertices in $V_1 \cup V_2$ are protected whereas vertices in $V_0$ are unprotected. The weight function $f$ is a \kaRD function if any subset $S\subseteq V_0$ of unprotected vertices are attacked, with $|S|\leq k$, there would be enough resources at the appropriate vertices in the open neighborhood of $S$ to defend the attack. 

The \textit{$k$-attack RD number} (\kaRD number) of a graph $G$, denoted $\gkaRD{G}$, is the minimum weight of all possible \kaRD functions on $G$.  If $\wf = \gkaRD{G}$ we will say that $f$ is a \textit{minimum weight} \kaRD \textit{function}. 

In this paper, we will primarily focus on \ptaRD, since the \oaRD  is equivalent to the original RD. We will also make generalizations to \kaRD when relevant. 

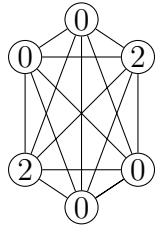
\begin{figure}[H]
\centering
\tikzstyle{vertex}=[circle,fill=black,inner sep=4.5pt]
\tikzstyle{white}=[circle,fill=white,inner sep=1pt]
\tikzstyle{green}=[circle,fill=green,inner sep=4.5pt]
\begin{tikzpicture}
 \node [draw, white] (1) at (0,0) {2};  
 \node [draw, white] (2) at (1.5,0) {0};
 \node [draw, white] (3) at (1.5,1.5) {2};
 \node [draw, white] (4) at (0,1.5) {0};
 \node [draw, white] (5) at (.75,2) {0};
 \node [draw, white] (6) at (.75,-.5) {0};
 \draw (1) -- (2) -- (3) -- (4) -- (5) -- (6) -- (1);
 \draw (1) -- (3) -- (5) -- (2) -- (4) -- (6) -- (2) -- (6) -- (3);
 \draw (1) -- (4);
 \draw (1) -- (5);
\end{tikzpicture}
\caption{Minimum Weight \taRD labeling of $K_6$}
\label{fig: K6}
\end{figure}

When studying dominating sets, complete graphs and star graphs are often used to demonstrate the efficiency of small dominating sets. While those graph classes also have small RD numbers, that is not the case for \ptaRD. For a complete graph of order $n$, denoted $K_n$, note that $\gamma(K_n) = 1$ and $\gamma_R(K_n) = 2$. Consider a $K_3$ and note that the minimum \taRD labeling is either a 0,1,2 or a 1,1,1 labeling of the three vertices. If we consider a $K_4$, there are more minimum \taRD labelings (0,2,0,2 or 1,1,1,1 for instance), all of which have a weight of 4. Notice that for any $K_n$ where $n \geq 4$, we have $\gtaRD{K_n} = 4$, since two vertices with label 2 are sufficient for a \taRD function. Note that if we only have one vertex with label 2, we can only have one vertex of label 0, and the rest of the vertices will have a label of 1. With two vertices of label 2, we have created a \taRD labeling.  A \taRD labeling of $K_6$ is demonstrated in Figure \ref{fig: K6} to illustrate this point.

This proves the following simple result: 

\begin{theorem}\label{2RDCompleteTheorem}
Let $K_n$ be a complete graph on $n$ vertices, then
$$\gtaRD{K_n} = \begin{cases} n, \text{ if } n < 4\\
4, \text{ if } n \geq 4.\\
\end{cases} $$
\end{theorem}

For positive integers $a$ and $b$, let $K_{a,b}$ denote the complete bipartite graphs with partitions of size $a$ and $b$. For the star graphs, note that $\gamma(K_{1,n}) = 1$ and $\gamma_R(K_{1,n}) = 2$. If we wanted to find a minimum weight \taRD labeling, we can begin by letting the center vertex have a label of 2. However, we will quickly note that while we can let one leaf vertex have a label of 0, all the other leaves must have a label of 1. With some exploration, we can easily show the following: 

\begin{theorem}\label{2RDStarTheorem}
Let $K_{1,n}$ be a star graph on $n+1$ vertices, then 

$$\gtaRD{K_{1,n}} = n+1.$$
\end{theorem}

As exemplified by the example above, it is not uncommon for a graph's \taRD number to be equal to the order of the graph. In the next section, we will introduce the idea of optimality to articulate this point.

\section{Optimality}
For most nontrivial graphs, the dominating and RD numbers are less than the order of the graph. However we could always give each vertex a label of 1.  This gives the following simple, but useful observation. 

\begin{prop} \label{prp: Gamma leq V(G)}
For any graph $G$, $\gkaRD{G} \leq |V(G)|$.
\end{prop}

\begin{figure}[H]
\centering
\tikzstyle{vertex}=[circle,fill=black,inner sep=4.5pt]
\tikzstyle{white}=[circle,fill=white,inner sep=1pt]
\tikzstyle{green}=[circle,fill=green,inner sep=4.5pt]
\centering
\begin{tikzpicture}
 \node [draw, white] (1) at (1.5,0) {2};  
 \node [draw, white] (2) at (3,1.5) {0};
 \node [draw, white] (3) at (2,1.5) {0};
 \node [draw, white] (4) at (1,1.5) {0};
 \node [draw, white] (5) at (0,1.5) {0};
 \node [draw, white] (6) at (1.5,3) {2};
 \node [draw, shape = circle, fill = black, minimum size = 0.04cm, inner sep=0pt] (7) at (1.4, 1.5) {};
 \node [draw, shape = circle, fill = black, minimum size = 0.04cm, inner sep=0pt] (7) at (1.5, 1.5) {};
 \node [draw, shape = circle, fill = black, minimum size = 0.04cm, inner sep=0pt] (7) at (1.6, 1.5) {};
 
 \draw (1) -- (2);
 \draw (1) -- (3);
 \draw (1) -- (4);
 \draw (1) -- (5);
 \draw (6) -- (2);
 \draw (6) -- (3);
 \draw (6) -- (4);
 \draw (6) -- (5);
\end{tikzpicture}
\caption{$K_{2,n}$ with the minimum weight \taRD labeling}
\label{fig: K1KnK1}
\end{figure}
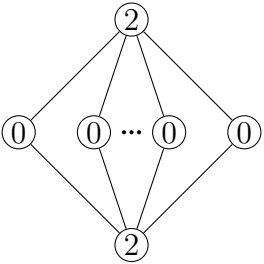

We want to differentiate between graphs for which $\gkaRD{G} < |V(G)|$ and $\gkaRD{G} = |V(G)|$. So, we define \textit{optimality} and \textit{optimal number} below.
\begin{definition}\label{def: optimality}
For a graph $G$, if $\gkaRD{G} = |V(G)|$, then $G$ is called \textit{sub-optimal}. If $\gkaRD{G} < |V(G)|$, then $G$ is called \textit{optimal}. If a graph is optimal then $|V(G)| - \gkaRD{G}$ is called the \textit{optimal number}.
\end{definition}
For the original domination function, the star graph is dominated by 1 single vertex. However, the same can not be said for \taRD. There is a similar graph class for \taRD that, like the star, can be 2-attack Roman Dominated by just two vertices: $K_{2,n}$ as shown in Figure \ref{fig: K1KnK1}. 

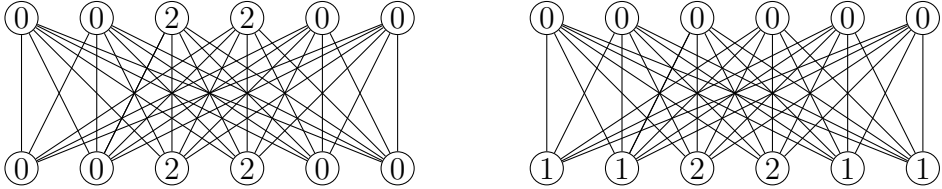
\begin{figure}[H]
\centering
\tikzstyle{vertex}=[circle,fill=black,inner sep=4.5pt]
\tikzstyle{white}=[circle,fill=white,inner sep=1pt]
\tikzstyle{green}=[circle,fill=green,inner sep=4.5pt]
\begin{subfigure}{0.40\textwidth}
\centering
\begin{tikzpicture}
 \node [draw, white] (1) at (3,0) {0};  
 \node [draw, white] (2) at (2,0) {0};
 \node [draw, white] (3) at (1,0) {2};
 \node [draw, white] (4) at (0,0) {2};
 \node [draw, white] (5) at (-1,0) {0};  
 \node [draw, white] (6) at (-2,0) {0};
 \node [draw, white] (9) at (3,2) {0};  
 \node [draw, white] (10) at (2,2) {0};
 \node [draw, white] (11) at (1,2) {2};
 \node [draw, white] (12) at (0,2) {2};
 \node [draw, white] (13) at (-1,2) {0};  
 \node [draw, white] (14) at (-2,2) {0};
 
 \draw (1) -- (9); \draw (1) -- (10); \draw (1) -- (11); \draw (1) -- (12); \draw (1) -- (13); \draw (1) -- (14);
 \draw (2) -- (9); \draw (2) -- (10); \draw (2) -- (11); \draw (2) -- (12); \draw (2) -- (13); \draw (2) -- (14);
 \draw (3) -- (9); \draw (3) -- (10); \draw (3) -- (11); \draw (3) -- (12); \draw (3) -- (13);\draw (3) -- (14);
 \draw (4) -- (9); \draw (4) -- (10);\draw (4) -- (11); \draw (4) -- (12); \draw (4) -- (13); \draw (4) -- (14); 
 \draw (5) -- (9); \draw (5) -- (10); \draw (5) -- (11); \draw (5) -- (12); \draw (5) -- (13); \draw (5) -- (14);
 \draw (6) -- (9); \draw (6) -- (10); \draw (6) -- (11); \draw (5) -- (12); \draw (6) -- (13); \draw (6) -- (14); 
\end{tikzpicture}
\end{subfigure}
\begin{subfigure}{0.40\textwidth}
\centering
\begin{tikzpicture}
 \node [draw, white] (1) at (3,0) {1};  
 \node [draw, white] (2) at (2,0) {1};
 \node [draw, white] (3) at (1,0) {2};
 \node [draw, white] (4) at (0,0) {2};
 \node [draw, white] (5) at (-1,0) {1};  
 \node [draw, white] (6) at (-2,0) {1};
 \node [draw, white] (9) at (3,2) {0};  
 \node [draw, white] (10) at (2,2) {0};
 \node [draw, white] (11) at (1,2) {0};
 \node [draw, white] (12) at (0,2) {0};
 \node [draw, white] (13) at (-1,2) {0};  
 \node [draw, white] (14) at (-2,2) {0};
 
 \draw (1) -- (9); \draw (1) -- (10); \draw (1) -- (11); \draw (1) -- (12); \draw (1) -- (13); \draw (1) -- (14);
 \draw (2) -- (9); \draw (2) -- (10); \draw (2) -- (11); \draw (2) -- (12); \draw (2) -- (13); \draw (2) -- (14);
 \draw (3) -- (9); \draw (3) -- (10); \draw (3) -- (11); \draw (3) -- (12); \draw (3) -- (13);\draw (3) -- (14);
 \draw (4) -- (9); \draw (4) -- (10);\draw (4) -- (11); \draw (4) -- (12); \draw (4) -- (13); \draw (4) -- (14); 
 \draw (5) -- (9); \draw (5) -- (10); \draw (5) -- (11); \draw (5) -- (12); \draw (5) -- (13); \draw (5) -- (14);
 \draw (6) -- (9); \draw (6) -- (10); \draw (6) -- (11); \draw (5) -- (12); \draw (6) -- (13); \draw (6) -- (14); 
\end{tikzpicture}
\end{subfigure}
\caption{A graph with multiple minimum \taRD labelings}
\label{figNonUnique2RD}
\end{figure}

Some graphs, such as the one in Figure \ref{fig: K1KnK1} have a unique minimum labeling, many graphs (see Figure \ref{figNonUnique2RD}) do not.

\begin{definition}
Let $\GtaRD{G}$ be the set of all minimum labelings on $G$:  $$\GtaRD{G} = \{f_G : \wfG = \gtaRD{G}\}.$$
\end{definition}

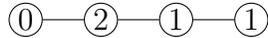
\begin{figure}[H]
\centering
\tikzstyle{vertex}=[circle,fill=black,inner sep=4.5pt]
\tikzstyle{white}=[circle,fill=white,inner sep=1pt]
\tikzstyle{green}=[circle,fill=green,inner sep=4.5pt]
\begin{tikzpicture}
 \node [draw, white] (1) at (0,2) {0};  
 \node [draw, white] (10) at (1,2) {2};
 \node [draw, white] (11) at (2,2) {1};
 \node [draw, white] (12) at (3,2) {1};
 \draw (1) -- (10) -- (11) -- (12);
\end{tikzpicture}
\caption{A graph with labeling 0-2-1-1.}
\label{shortpath}
\end{figure}

Notice that both labelings in Figure \ref{figNonUnique2RD} are minimum and thus would be elements of $\GtaRD{G}$. Throughout we will assume that $\lG:V(G) \rightarrow \{0,1,2\}$ will denote a labeling in $\gtaRD{G}$. So $w(\lG) = \gtaRD{G}$. When discussing labels of sequential vertices in a graph, we will use dashes to delineate the subpath of vertices we wish to highlight. For instance, the graph in Figure \ref{shortpath} has a subgraph labeled 0-2-1 and a subgraph labeled 1-1-2-0, but they are not disjoint subpaths. 

Here we use exterior public and private neightbors to discuss how $V_2(f)$ interacts with $V_0(f)$.

\begin{definition}\label{def: V0V2neighbors}
Given a graph $G$ and labeling $f_G$, if a vertex $v \in V_0(f_G)$ is adjacent to exactly one vertex in $V_2(f_G)$, then we say $v \in epn(V_2(f_G))$. If a vertex $u \in V_0(f_G)$ is adjacent to more than one vertex in $V_2(f_G)$, then we say $u \in \overline{ep}n(V_2(f_G))$. 
\end{definition}

For emphasis, we are only referring to vertices with label 0 for these neighborhoods. A vertex with label 1 will never be considered in the exterior public or private neighborhood of $V_2$. 

\subsection{What Graphs are Optimal?}
In this section, we will discuss what properties a graph must have to be optimal, and how to find the \taRD number of any graph. The following proposition is intuitive but important.

\begin{prop}\label{prp: V_0 > V_2}
A graph $G$ is optimal if and only if there is a \taRD labeling $\ell \in \GtaRD{G}$ so that $|V_2(\ell)| < |V_0(\ell)|$.
\end{prop}

Note that if there exists a \taRD labeling $\ell \in \GtaRD{G}$ with $|V_2(\ell)|<|V_0(\ell)|$ then $|V_0(\ell)| - |V_2(\ell)| = |V(G)| - \gtaRD{G}$.  Thus any other \taRD labeling $\ell' \in \GtaRD{G}$ will have $|V_2(\ell')|<|V_0(\ell')|$. We use this idea to help prove the following theorem.

\begin{theorem}\label{P5Subgraph}
A graph $G$ is optimal if and only if for every labeling $\ell_G \in \GtaRD{G}$, there exists a $P_5$ subgraph in $G$ labeled 0-2-0-2-0 under $\ell_G$.
\end{theorem}

\begin{proof}

Part 1: 
Assume $G$ has a $P_5$ subgraph, $abcde$, which is labeled 0-2-0-2-0 under $\ell$. Let $f(a) = f(c) = f(e) = 0$. Let $f(b) = f(d) = 2$. Let $f(v) = 1$ for all other $v \in V(G)$. This is a valid \taRD labeling, where $\w{\ell} \leq \wfG <  |V(G)|$. Since the weight of any valid \taRD labeling is greater than or equal to $\gtaRD{G}$, it follows that $\gtaRD{G} < |V(G)|$ and thus $G$ is optimal. 

Part 2: 
Let $\ell \in \GtaRD{G}$ and assume $G$ is optimal.  By Proposition \ref{prp: V_0 > V_2}, we know $|V_0|>|V_2|$. Since $V_0$ is non-empty, $V_2$ must also be non-empty for the labeling to be valid. If every $v \in V_2$ is only adjacent to one vertex with label 0, then $|V_0| \leq |V_2|$ which is a contradiction. So, there must be a $b \in V_2$ where $b$ is adjacent to at least two vertices with label 0 (call them $a$ and $c$). Since $\ell$ is a valid labeling, either $a$ or $c$ ($wlog$ assume it is $c$) must be adjacent to another vertex with label 2 (call this vertex $d$). If $d$ is adjacent to another vertex with labeled 0, we have found a $P_5$ with a 0-2-0-2-0 labeling and we are done. If not, relabel $d$ and $c$ to both have a label of 1 (as shown in Figure \ref{2-0 to 1-1}). Call this labeling $\ell^1$. Note, $\w{\ell} = \w{\ell^1}$ so $\ell^1 \in \GtaRD{G}$, which tells us that $|V_0(\ell^1)|=|V_0(\ell)|-1$ and $|V_2(\ell^1)|= |V_2(\ell)|-1$ and so $|V_0(\ell^1)|-|V_2(\ell^1)| = |V_0(\ell)|-|V_2(\ell)|>0$. 

Note that $\ell$ and $\ell^1$ are equal on all vertices in $V_0(\ell^1) \cup V_2(\ell^1)$. So, if there is a 0-2-0-2-0 $P_5$ subgraph in $G$ under $\ell^1$, then the same $P_5$ subgraph and labeling will exists under $\ell$. Find another vertex, $b^1 \in V_2(\ell^1)$ that has two adjacent vertices with label 0 (call them $a^1$ and $c^1$) and repeat this process. If we reach a point where every vertex in $V_2(\ell^n)$ is adjacent to at most one vertex with label 0, as stated previously, there is no way to have $|V_0(\ell)| > |V_2(\ell)|$, which violates our assumption that $G$ is optimal.

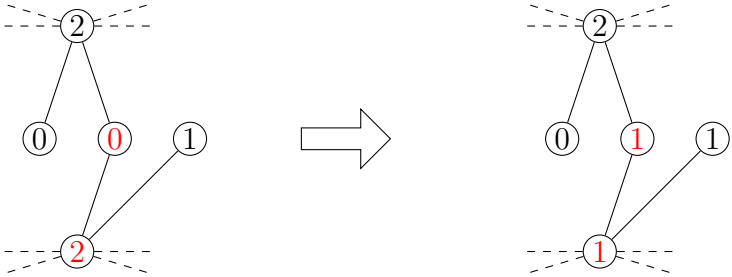
\begin{figure}[H]
\centering
\tikzstyle{vertex}=[circle,fill=black,inner sep=4.5pt]
\tikzstyle{white}=[circle,fill=white,inner sep=1pt]
\tikzstyle{green}=[circle,fill=green,inner sep=4.5pt]
\begin{subfigure}{0.40\textwidth}
\centering
\begin{tikzpicture}
 \node [draw, white] (1) at (1.5,0) {\rc{2}};  
 \node [draw, white] (2) at (3,1.5) {1};
 \node [draw, white] (3) at (2,1.5) {\rc{0}};
 \node [draw, white] (4) at (1,1.5) {0};
 \node [draw, white] (6) at (1.5,3) {2};
 \node [My Arrow Style] at (6,1.5) {};
 \draw (1) -- (2);
 \draw (1) -- (3);
 \draw (1) -- (.5,0) [dashed];
  \draw (1) -- (.5,-.3) [dashed];
  \draw (1) -- (2.5,-.3) [dashed];
  \draw (1) -- (2.5,0) [dashed];
   \draw (6) -- (.5,3) [dashed];
  \draw (6) -- (.5,3.3) [dashed];
  \draw (6) -- (2.5,3.3) [dashed];
  \draw (6) -- (2.5,3) [dashed];
 \draw (6) -- (3);
 \draw (6) -- (4);
\end{tikzpicture}
\end{subfigure}
\begin{subfigure}{0.4\textwidth}
\centering
\begin{tikzpicture}
 \node [draw, white] (1) at (1.5,0) {\rc{1}};  
 \node [draw, white] (2) at (3,1.5) {1};
 \node [draw, white] (3) at (2,1.5) {\rc{1}};
 \node [draw, white] (4) at (1,1.5) {0};
 \node [draw, white] (6) at (1.5,3) {2};
 \draw (1) -- (2);
 \draw (1) -- (3);
 \draw (1) -- (.5,0) [dashed];
  \draw (1) -- (.5,-.3) [dashed];
  \draw (1) -- (2.5,-.3) [dashed];
  \draw (1) -- (2.5,0) [dashed];
  \draw (6) -- (.5,3) [dashed];
  \draw (6) -- (.5,3.3) [dashed];
  \draw (6) -- (2.5,3.3) [dashed];
  \draw (6) -- (2.5,3) [dashed];
 \draw (6) -- (3);
 \draw (6) -- (4);
\end{tikzpicture}
\end{subfigure}
\caption{Turning a 2-0 pair into a 1-1 pair}
\label{2-0 to 1-1}
\end{figure}
\end{proof}

Theorem \ref{P5Subgraph} is incredibly useful, as it is a quick way not only to tell whether a graph is optimal or not, but gives us a start for what an optimal \taRD labeling would look like. In the next section, we will consider how to find the minimum \taRD labeling of any graph.

\subsection{Finding the Minimum labeling}
In order to find a minimum \taRD labeling of a graph $G$, we must first introduce the concept of end-coupled center disjoint $P_5$ subgraphs. 

\begin{definition}\label{ECCDP5}
Let $S$ be a set of $P_5$ subgraphs of $G$ (not necessarily disjoint). For $S$ to be a valid  \textit{end-coupled center disjoint (eccd) $P_5$ set}, we must have:
\begin{itemize}
    \item End Coupled: If $abcde$ is some $P_5$ in $S$ and $a$ or $b$ is a vertex in some other $P'_5$ in $S$, then $ab \in E(P'_5)$, where $a$ is a leaf. Analogously, if $e$ or $d$ is a vertex in some other $P''_5$ in $S$, then $ed \in E(P''_5)$, where $e$ is a leaf.
    \item Center Disjoint: If $abcde$ is some $P_5$ in $S$ then no other $P_5$ in $S$ may use the vertex $c$.
\end{itemize}
\end{definition}
So, if every $P_5$ in $S$ is completely disjoint, $S$ is a valid set of eccd $P_5$ set. If there is overlap (vertices appear in multiple subgraphs in $S$), then the overlap must be done ``correctly" in the manner described. Because we will only be concerned with $P_5$ subgraphs, we will often say a set $S$ is a valid eccd set. 
\begin{figure}[H]
\centering
\tikzstyle{vertex}=[circle,fill=black,inner sep=2pt]
\centering
\begin{tikzpicture}
 \node [draw, vertex, label = left:$v_1$] (1) at (-2, 0) {};  
 \node [draw, vertex, label = above:$v_2$] (2) at (-.5,0) {};
 \node [draw, vertex, label = below:$v_8$] (3) at (1,-1) {};
  \node [draw, vertex, label = above:$v_4$] (4) at (2.5,0) {};
  \node [draw, vertex, label = right:$v_5$] (5) at (4,0) {};  
   \node [draw, vertex, label = left:$v_6$] (6) at (-2, -1) {};  
 \node [draw, vertex, label = below:$v_7$] (7) at (-.5,-1) {};
  \node [draw, vertex, label = below:$v_9$] (8) at (2.5,-1) {};
  \node [draw, vertex, label = right:$v_{10}$] (9) at (4,-1) {};
  \node [draw, vertex, label = above:$v_{3}$] (10) at (1,0) {};

\draw (1) -- (2) -- (3) -- (4) -- (5);
\draw (6) -- (7) -- (3) -- (8) -- (9);
\draw (2) -- (10) -- (4);
\draw (5) -- (3);
\end{tikzpicture}
\caption{A graph with several possible eccd sets}
\label{fig:eccdExample}
\end{figure}
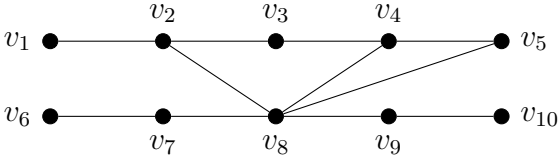
Consider the graph in Figure \ref{fig:eccdExample}. Obviously, $v_1v_2v_3v_4v_5$ and $v_6v_7v_8v_9v_{10}$ would make up a valid eccd set as the two $P_5's$ are disjoint. However, $v_1v_2v_8v_4v_5$ and  $v_6v_7v_8v_9v_{10}$ would \textbf{not} be a valid eccd set as it violates the center disjoint criteria. Additionally, notice that $v_1v_2v_3v_4v_5$ and $v_6v_7v_8v_4v_5$ is a valid eccd set as the $P_5's$ overlap on $v_4v_5$, but that end coupling is held through both subgraphs. Alternatively, the subgraphs $v_1v_2v_3v_4v_5$ and $v_{10}v_9v_8v_5v_{4}$ would \textbf{not} be a valid eccd set as $v_4v_5$ are coupled in both, but appear in the wrong order ($v_5$ is a leaf in one subgraph and $v_4$ is a leaf in the other subgraph). 

For the graph in Figure \ref{fig:eccdExample}, there are several valid eccd sets, the largest sets having size 2. We can convince ourselves that no eccd set of size 3.  

Before implementing the use of eccd sets, we must first consider the following propositions and lemma on how the sets $V_0$ and $V_2$ interact.

\begin{prop}\label{prop:edgeremoval}
    For any graph $G$ with edge $e$, $$\gtaRD{G-e} \geq \gtaRD{G}.$$
\end{prop}

\begin{proof}
    If $f \in \GtaRD{G-e}$, then $f$ is a valid \taRD function on $G$, but may not be minimal.  The result follows. 
\end{proof}

\begin{prop}\label{prop: 2adj0}
Given a graph $G$ with $\ell_G \in \GtaRD{G}$, every vertex in $V_2$ is adjacent to at least one vertex in $V_0$.
\end{prop}

\begin{proof}

If not, we can relabel any vertex in $V_2$ which is not adjacent to a vertex in $V_0$ with the label 1, thus reducing the total weight. 
\end{proof}

Proposition \ref{prop: 2adj0} leads to the following lemma.


\begin{lemma}\label{lma: EveryEPN}
Given a graph $G$ with labeling $\ell_G \in \GtaRD{G}$, there exists a matching (independent edge set), $M$, so that so that every $v \in V_2$ is matched with a $v' \in  V_0 \cap epn(V_2)$. In other words, there exists a bijection between $V_2$ and some $V'_0 \subseteq V_0$.
\end{lemma}

\begin{proof}
Consider a graph $G$ with $\ell_G \in \GtaRD{G}$. Consider a matching, $M_0$, in $G$ so that every vertex $v\in V_2$ is matched with a vertex in $epn(v, V_2)$ as long as $|epn(v,V_2) \cap V_0| \geq 1$. Now, label all the remaining vertices in $V_2$ as $t_1, t_2, \ldots, t_k$.  Note that each $t_i$ has no exterior private neighbor in $V_0$ but must be adjacent to at least one vertex in $V_0$ by Proposition \ref{prop: 2adj0}. Pick a $v_1 \in V_0$ adjacent to $t_1$. Note that $v_1$ could not be in $M_0$ since $v_1$ is not an exterior private neighbor of $V_2$.  Consider the spanning subgraph $G_1$ which is $G$ with all the edges from $v_1$ to any vertex in $V_2 \backslash t_1$ removed. Notice in $G_1$, $v_1$ is the $epn$ of $t_1$. The induced labeling, $\ell_{G_1}$ on $G_1$ is valid by construction and has the same weight as $\ell_G$ thus is still minimum by Proposition \ref{prop:edgeremoval}. Add the matching $v_1t_1$ to $M_0$ to form $M_1$.  

Note Proposition \ref{prop: 2adj0} still holds for $G_1$, so we can repeat this process for every vertex $t_i$ for $2 \leq i \leq k$. We will end with a subgraph $G_k$ with a minimum \taRD labeling $\ell_{G_k}$ where every vertex in $V_2$ has an $epn$ in $G_k$ producing a matching $M_k$ with the desired condition. 
\end{proof}

Now we can proceed to find the optimal number.  Given a graph $G$, let $\mathbb{S}= \{S: S \text{ is a eccd set of } G \}$.
\begin{theorem}\label{thm: MAINECCD}
The optimal number (recall Definition \ref{def: optimality}) of a graph $G$ is equal to the maximum number of distinct end-coupled center disjoint $P_5$ subgraphs in $G$:
 $$ |V(G)| - \gtaRD{G} =\max\{|S|:S \in \mathbb{S}\}.$$
\end{theorem}

Before we prove this theorem, it is helpful to note that a labeling constructed from a set of eccd $P_5$ subgraphs is valid.

\begin{lemma} \label{lma: eccdValidLabeling}
 Given a graph $G$ and a set of $S$ of eccd subgraphs in $G$, define a labeling $f_G$ where we label every eccd subgraph in $S$ as $0-2-0-2-0$ (some vertices may receive a label more than once, but each time it will be the same value by construction). Give every vertex in $G \backslash S$ a label of 1. Then $f_G$ is a valid labeling.
\end{lemma}

\begin{proof}

Let $f_G$ be a labeling on $G$ as described. We can break down the definition of \taRD into two parts: (1) every 0 must be adjacent to a 2, and (2) every pair of two 0's must be adjacent to at least two 2's. 

Notice that if you label eccd $P_5$'s in the manner described, and those are the only vertices in $V_0(f)$, then every 0 will be adjacent to at least one 2. Now, assume for the sake of contradiction, there are two vertices $a,b \in V_0(f)$ which are both only adjacent to $c \in V_2(f)$. Neither $a$ or $b$ can be a center vertex of an eccd $P_5$.  If $a$ and $b$ were both pendent vertices of two $P_5$'s, then $c$ was incorrectly coupled with two vertices, which contradicts how eccd sets are defined.
\end{proof}

Now, to prove Theorem \ref{thm: MAINECCD}, we will show that the maximum number of distinct eccd $P_5$ subgraphs is equal to the difference $|V_0| - |V_2|$ of $\lG$. For ease of notation, let $\max\{|S|:S \in \mathbb{S}\} = \mathbf{E}(G)$.
\begin{proof} (Proof of Theorem \ref{thm: MAINECCD})

Given a graph $G$ with labeling $\lG \in \GtaRD{G}$.  Construct a spanning subgraph $G'$ with $\lGG = \lG$ for all $v \in V(G')$, such that every $v \in V_2(\lGG)$ has an $epn$. Let $V'_0(\lGG)$ denote the set of $epn$'s of $V_2(\lGG)$ noting that $|V_0'(\lGG)| = |V_2(\lGG)|$. This construction is possible by Lemma \ref{lma: EveryEPN}. Now let $V''_0(\lGG)  = V_0(\lGG) \backslash V'_0(\lGG)$. By construction, every vertex in $V''_0(\lGG)$ is a public neighbor of $V_2(\lGG)$. Recall that the optimal number of $G$ is $|V_0(\ell_{G})| - |V_2(\ell_{G})|$ and  $|V_0(\ell_{G})| - |V_2(\ell_{G})| =  |V''_0(\lGG)| + |V'_0(\lGG)| - |V_2(\lGG)|$.  Therefore, the optimal number of $G$ is $|V''_0(\lGG)|$.

So, we now will show that $|V''_0(\lGG)|$ is equal to the number of eccd $P_5$ subgraphs in $G$. First, note we can take the $\max\{|S|:S \in \mathbb{S}\}$ and use said $S$ to create a labeling $f_{G'}$ in the manner described in Lemma \ref{lma: eccdValidLabeling}. We know this labeling is valid, and can deduce that $|V(G')| - w(f_{G'}) = \mathbf{E}(G)$. Therefore if $\mathbf{E}(G)$ was greater than $V''_0(\lGG)$, then $w(f_{G'}) < w(\lGG)$, which is a contradiction as $\lGG$ is assumed to be minimum. This implies $|V''_0(\lGG)| \geq \mathbf{E}(G)$.

Note that every vertex in $V''_0(\lGG)$ is the center of a $P_5$ since each vertex in $V''_0(\lGG)$ must be adjacent to at least two vertices in $V_2(\lGG)$ and each vertex in $V_2(\lGG)$ is coupled with a vertex in $V'_0(\lGG)$. Therefore, each of these $P_5$'s are center disjoint.  In addition, the $P_5$'s are end-coupled by construction. This implies $|V''_0(\lGG)| \leq \mathbf{E}(G)$. Therefore, $|V''_0(\lGG)| = \mathbf{E}(G)$.
\end{proof}

We have shown that finding the maximum number of eccd $P_5$ subgraphs tells us the optimal number of a graph, but in fact, it does more than that. Lemma \ref{lma: eccdValidLabeling} tells us that we can construct a valid labeling where every eccd $P_5$ subgraph is labeled $0 - 2 - 0 -2 -0$ and all other vertices receiving a label of 1. Not only will a labeling of this manner be valid, but it will also be minimal.

\begin{corr}\label{corr: optimallabeling}
The maximum set of eccd $P_5$ subgraphs of a graph $G$ gives us a minimum weight \taRD labeling of $G$, $\ell^*_G$, where $\w{\ell^*_G} = \gtaRD{G}$.
\end{corr}

\begin{proof}

Given a graph $G$ we can use Lemma \ref{lma: eccdValidLabeling} to create a labeling $\ell^*_G$. Let $\ell^*_G$ be a labeling where we take the largest set $S$ of eccd $P_5$'s in $G$ and label each $0-2-0-2-0$ with all other vertices receiving a label of 1. Lemma \ref{lma: eccdValidLabeling} ensures this labeling is valid. As discussed in the proof of Theorem \ref{thm: MAINECCD}, we know $|V(\ell^*_G)| - w(\ell^*_G) = |S|$. (Recall $|V(\ell^*_G)| - w(\ell^*_G) = |V_0(\ell^*_G)| - |V_2(\ell^*_G)|$ and every eccd adds one to the latter difference.) Since $S$ is maximum, it follows that 
\begin{align*}
    |V(G)| - \gtaRD{G}  &= |V(\ell^*_G)| - w(\ell^*_G) \\ \gtaRD{G} &= w(\ell^*_G).
\end{align*}
\end{proof}

Often, the minimum weight labeling is not unique, so this provides one valid minimum labeling to any graph. 


\section{Density of 2-attack Roman Domination}

In this section, we will discuss the \taRD density of a graph which is another way to describe the nature of how the \taRD function behaves on a graph. This could easily extend the idea of \taRD number to infinite graphs. 

Following Theorem \ref{thm: MAINECCD}, we can very quickly find that the optimal number for both paths and cycles is $\floor{\frac{n}{5}}$. For example, consider $C_{17}$. We can label three consecutive disjoint $P_5$ subgraphs $0-2-0-2-0$ and then place a 1 on the remaining two vertices. Using this logic and Corollary \ref{corr: optimallabeling}, we find that
\begin{align*} 
   \gtaRD{C_n} = \gtaRD{P_n} &= 4\floor{\frac{n}{5}} + \left(n \mod \text{ } 5\right)\\
   &= 4\floor{\frac{n}{5}} + \left(n - 5\floor{\frac{n}{5} }\right)\\
   &= n - \floor{\frac{n}{5}}.
\end{align*}
Now we ask the question, what is the \taRD number relative to the order of the graph?  This helps us characterize how efficient or inefficient certain graphs are.

\begin{definition}\label{defn: Density}
The \taRD \textit{density} of a graph $G$ is 
$$ \mathbf{D}(G) = \frac{\gtaRD{G}}{|V(G)|}. $$
\end{definition}
Intuitively, a graph like the one in Figure \ref{fig: K1KnK1} will have a low density, as we can have a graph with very large order and a label sum of 4. Ultimately, a low density is desirable as that is indicative of how efficient a minimum labeling can be. But, how low can the \taRD dominating density be?

\begin{theorem}\label{thm: densitybound}
The \taRD dominating density of a graph $G$ is bounded below by 
$$ \mathbf{D}(G) \geq \frac{4}{\Delta(G) + 3}. $$
\end{theorem}
In order to prove this theorem, we first need two short lemmas:

\begin{lemma}\label{lma: G'labelingMustBeG}
 Given a graph $G$ with minimum labeling $\ell_G$, let $G'$ be an induced subgraph $G$ in which $V(G') = V(G) \backslash V_1(\ell_G)$. Then, the induced labeling from $G$ is a minimum labeling of $G'$.
\end{lemma}

\begin{proof}

If there were a smaller weight minimum labeling of $G'$ then we can extend that labeling to $G$ by labeling all vertices in $V_1(\ell_G)$ as 1.  This is a new labeling on $G$ whose total weight is less than that of $\ell(G)$ which is a contradiction. 
\end{proof}

Using this idea, we are able to prove another property of our induced $G'$ subgraph.  The following lemma shows that if we have a graph with a minimum labeling and remove all the vertices which have a label 1, then the density can not increase.  This seems intuitive since the maximum value for the density of any graph is 1 (by Proposition \ref{prp: Gamma leq V(G)}). But this is often not the minimum labeling, which would only decrease the density, since $\gtaRD{G} \leq |V(G)|$.

\begin{lemma}\label{lma: densitySubgraph}
 Given a graph $G$ with minimum labeling $\ell_G$, let $G'$ be the induced subgraph on $V(G) \backslash V_1(\ell_G)$. Then, $\mathbf{D}(G') \leq  \mathbf{D}(G)$.
\end{lemma}

\begin{proof}

Consider the case where $G$ is not optimal. Thus $\mathbf{D}(G) = 1$. Since $\gtaRD{H} \leq |V(H)|$ for any graph $H$, we know $\mathbf{D}(H)\leq 1$ for any graph $H$.  Thus $\mathbf{D}(G')\leq 1 = \mathbf{D}(G)$.

Assume that $G$ is optimal, thus $\mathbf{D}(G) < 1$. Since $G'$ is obtained by removing every vertex in $G$ that is labeled 1 in $\ell_G$, 
\begin{align*}
   \mathbf{D}(G') &=  \frac{\gtaRD{G} - |V_1(\ell_G)|}{|V(G)| - |V_1(\ell_G)|}\\
   & \leq \frac{\gtaRD{G}}{|V(G)|}\\
   & = \mathbf{D}(G),
\end{align*}
since $\frac{a-b}{n-b} \leq \frac{a}{n}$ for $0 \leq b \leq a < n$.  
\end{proof}

So, Lemmas \ref{lma: G'labelingMustBeG} and \ref{lma: densitySubgraph} inform us that we can remove the vertices with label 1 in some $\lG$, and the resulting labeling is valid, minimal, and will have not have a greater density than the original graph. Using this, we are now ready to prove Theorem \ref{thm: densitybound}. Recall that $\Delta(G)$ is the maximum degree of $G$. For ease of notation in the proof below, let $k = \Delta(G)$.

\begin{proof} (Proof of Theorem \ref{thm: densitybound})

Let $G$ be a graph with $\ell_G \in \gtaRD{G}$. Let $k = \Delta(G)$. Let $G'$ be the induced subgraph in which $V(G') = V(G)\backslash V_1(\ell_G)$ and let $\ell_{G'}$ be a labeling on $G'$ induced by $\ell_G$. By Lemma \ref{lma: densitySubgraph}, showing that $\mathbf{D}(G') \geq \frac{4}{k + 3}$ is sufficient to prove that $\mathbf{D}(G) \geq \frac{4}{k+3}$.  

Let $V_2=V_2(\ell_{G'})$ and $V_0=V_0(\ell_{G'})$. Define $A_2=\{v \in V_2: v \text{ has an epn in } V_0\}$.  Let $A_0$ be the set of epn for those vertices in $A_2$.  A vertex in $A_2$ can not have two epn in $A_0$ since $\ell_{G'}$ is a valid \taRD labeling, therefore $|A_0| = |A_2|$. Let $B_2 = V_2/A_2$ and let $B_0 = V_0/A_0$. 

Note that
\begin{align*}
     \mathbf{D}(G') = \frac{\gtaRD{G'}}{|V(G')|} = \frac{2|V_2|}{|V_2| + |V_0|}.
\end{align*}

Therefore, if
$$ |V_0| \leq \frac{|V_2|}{2}(k-1) + |V_2|,$$ then one can deduce
$$\frac{4}{k+3} \leq \mathbf{D}(G'). $$

Every vertex in $A_0$ is incident to exactly one vertex in $V_2$.  Every vertex in $B_0$ is incident to at least two vertices in $V_2$. Therefore, if $x$ is the number of edges from $V_2$ to $B_0$ then 

$$2|B_0| \leq x \leq |A_0|(k-1) + k|B_2|.$$
This implies
\begin{align*}
    |V_0| &\leq \frac{|A_0|}{2}(k-1) + \frac{k|B_2|}{2} + |A_0| \\
        &=\frac{k-1}{2}\left(|A_0|+|B_2|\right) + \frac{|B_2|}{2} + |A_0|\\
        &\leq \frac{k-1}{2} |V_2| + |V_2|
\end{align*}
which completes the proof. 
\end{proof}

Recall that we assumed that $|V_1(\lG)| \geq 0$ in the proof for Lemma \ref{lma: densitySubgraph}. If we assume $|V_1(\lG)| > 0$, that is, there are vertices with label 1 in the minimum labeling, then we will arrive at a strict inequality $\frac{a}{b} > \frac{a-n}{b-n}$ at the end of the proof for Lemma \ref{lma: densitySubgraph}. Thus, the following corollary emerges.

\begin{corr}
If a graph $G$ is optimal, and there exists a minimum labeling of $G$, $\lG$ such that $|V_1(\lG)| > 0$, then $\mathbf{D}(G) > \frac{4}{\Delta(G) + 3}.$
\end{corr}

\begin{proof}

Follows from the proof of Lemma \ref{lma: densitySubgraph}.
\end{proof}

\subsection{Maximum Degree Bound}\label{sec: MaxDegBound}
It can be easily shown that given a graph $G$ with order $n$, then $\gamma(G) \geq \frac{n}{\Delta(G) + 1}$. This is a rather simple concept to understand because every vertex that is in the dominating set is adjacent to at most $\Delta(G)$ vertices, which means that a vertex can dominate at most $\Delta(G) + 1$ vertices, including itself.  Therefore, the minimum number of vertices we need in a dominating set is the order of the graph divided by $\Delta(G) + 1$.  
 
Cockayne et. al \cite{Cockayne} proved that $\gR{G} \geq \frac{2n}{\Delta(G) + 1}$. Now, if we re-write the lower bound for density as shown in Theorem \ref{thm: densitybound}, we get that $\gtaRD{G} \geq \frac{4n}{\Delta(G) + 3}$. Observing these three lower bounds in sequence makes the progression clear, and said progression makes intuitive sense when we unpack how each dominating variant behaves differently than the last. Under Roman Domination, a vertex with label 2 can fully dominate vertices adjacent to it. Whereas in \ptaRD, a vertex with label 2 can only partially dominate vertices adjacent to it, but it needs the help of another vertex of label 2 to make the labeling valid. So, we see how, much like the jump from domination to Roman Domination, we need almost \textit{double} the amount of dominating vertices to have a valid \taRD function. Alas, the extra two vertices in the denominator, as seen more clearly by rewriting the bound for $\gtaRD{G}$ as $\gtaRD{G} \geq \frac{4n}{\Delta(G) + 1 + 2}$ come from the fact that each 2 can fully dominate at most one vertex (their $epn$).
 
\subsection{Exterior Private Neighbors and Minimum Weight Labeling}
We have discussed throughout that we typically strive to label a graph such that vertices with label 2 have exterior private neighbors, as it extends the reach that dominating vertices have on the graph. However, this is not always achievable with certain graphs; for example, the graphs in Figure \ref{fig: K1KnK1} and Figure \ref{figNonUnique2RD} both have minimum weight \taRD labelings with no $epn$'s. So, we cannot argue that a minimum weight labeling implies 2's have $epn$'s. But one may want to argue that the converse is true: If every vertex in some $f$ on $G$ with label 2 has an $epn$, then $f$ is a minimum weight labeling. However, we can see that this is not the case in the following example.
\begin{figure}[H]
\centering
\tikzstyle{vertex}=[circle,fill=black,inner sep=4.5pt]
\tikzstyle{white}=[circle,fill=white,inner sep=1pt]
\tikzstyle{green}=[circle,fill=green,inner sep=4.5pt]
\centering
\begin{tikzpicture}
 \node [draw, white] (1) at (-2,2) {2};  
 \node [draw, white] (2) at (-2,0) {2};
 \node [draw, white] (3) at (-2,-2) {2};
 
 \node [draw, white] (4) at (0,1) {0};  
 \node [draw, white] (5) at (0,-1) {0};
 
 \node [draw, white] (6) at (2,2) {\rc{0}};  
 \node [draw, white] (7) at (2,0) {\rc{0}};
 \node [draw, white] (8) at (2,-2) {\rc{0}};
 
 \draw (1) -- (6);
 \draw (2) -- (7);
 \draw (3) -- (8);
 
 \draw (4) -- (1);
  \draw (4) -- (2);
   \draw (4) -- (3);
    \draw (4) -- (6);
     \draw (4) -- (7);
      \draw (4) -- (8);
      
\draw (5) -- (1);
  \draw (5) -- (2);
   \draw (5) -- (3);
    \draw (5) -- (6);
     \draw (5) -- (7);
      \draw (5) -- (8);
\end{tikzpicture}
\caption{A labeling on $G$ where every $v \in V_2$ has an $epn$}
\label{fig: Allepns}
\end{figure}
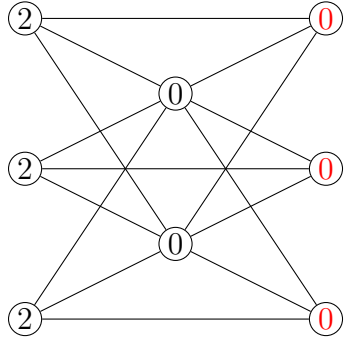
In the graph in Figure \ref{fig: Allepns}, all three vertices with label 2 have an exterior private neighbor (highlighted in red). However, this labeling is not minimum.
\begin{figure}[H]
\centering
\tikzstyle{vertex}=[circle,fill=black,inner sep=4.5pt]
\tikzstyle{white}=[circle,fill=white,inner sep=1pt]
\tikzstyle{green}=[circle,fill=green,inner sep=4.5pt]
\centering
\begin{tikzpicture}
 \node [draw, white] (1) at (-2,2) {0};  
 \node [draw, white] (2) at (-2,0) {0};
 \node [draw, white] (3) at (-2,-2) {0};
 
 \node [draw, white] (4) at (0,1) {2};  
 \node [draw, white] (5) at (0,-1) {2};
 
 \node [draw, white] (6) at (2,2) {0};  
 \node [draw, white] (7) at (2,0) {0};
 \node [draw, white] (8) at (2,-2) {0};
 
 \draw (1) -- (6);
 \draw (2) -- (7);
 \draw (3) -- (8);
 
 \draw (4) -- (1);
  \draw (4) -- (2);
   \draw (4) -- (3);
    \draw (4) -- (6);
     \draw (4) -- (7);
      \draw (4) -- (8);
      
\draw (5) -- (1);
  \draw (5) -- (2);
   \draw (5) -- (3);
    \draw (5) -- (6);
     \draw (5) -- (7);
      \draw (5) -- (8);
\end{tikzpicture}
\caption{A minimum labeling of the graph in Figure \ref{fig: Allepns} where no $v \in V_2$ has an $epn$}
\label{fig: Noepns}
\end{figure}
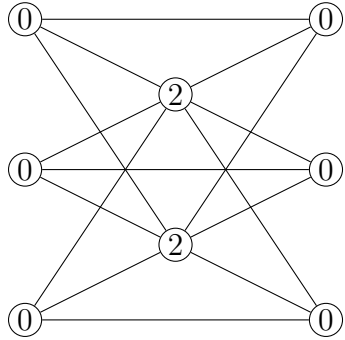
Note the labeling of the graph in Figure \ref{fig: Noepns} has a smaller weight than the labeling of the same graph in Figure \ref{fig: Allepns} but has no vertices in $V_2$ with exterior private neighbors. So, our brief analysis tells us that although exterior private neighbors can be helpful, and in cases result in a minimum labeling, there is no immediate relationship between maximizing the number of 2's which have $epn$'s and minimizing the weight of the \taRD function on a graph.

\bibliographystyle{abbrv} 
\bibliography{BibFile}

\end{document}